\documentclass[11pt,reqno]{article}
\usepackage{latexsym, amsmath, amssymb, a4, epsfig}
\usepackage{graphicx}

\newtheorem{theorem}{Theorem}[section]

\newcommand{\nm}{\noalign{\smallskip}}
\newcommand{\qed}{ $\Box$}
\newcommand{\ds}{\displaystyle}
\newcommand{\pf}{\noindent {\sl Proof}. \ }
\newcommand{\p}{\partial}
\newcommand{\pd}[2]{\frac {\p #1}{\p #2}}
\newcommand{\eqnref}[1]{(\ref {#1})}

\newcommand{\Rbb}{\mathbb{R}}

\newcommand{\la}{\langle}
\newcommand{\ra}{\rangle}

\newcommand{\Acal}{\mathcal{A}}

\newcommand{\Ccal}{\mathcal{C}}

\newcommand{\Hcal}{\mathcal{H}}

\newcommand{\Kcal}{\mathcal{K}}

\newcommand{\Scal}{\mathcal{S}}

\newcommand{\Ga}{\alpha}

\newcommand{\Gd}{\delta}
\newcommand{\Ge}{\epsilon}

\newcommand{\Gvf}{\varphi}

\newcommand{\Gl}{\lambda}

\newcommand{\Gs}{\sigma}

\newcommand{\Go}{\omega}

\newcommand{\GD}{\Delta}

\newcommand{\GG}{\Gamma}

\newcommand{\GO}{\Omega}

\newcommand{\beq}{\begin{equation}}
\newcommand{\eeq}{\end{equation}}

\numberwithin{equation}{section}
\numberwithin{figure}{section}

\begin{document}
\title{Analysis of plasmon resonance on smooth domains using spectral properties of the Neumann-Poincar\'e operator\thanks{\footnotesize This work is supported by the Korean Ministry of Education, Sciences and Technology through NRF grants No. 2010-0017532.}}

\author{Kazunori Ando\thanks{Department of Mathematics, Inha University, Incheon
402-751, S. Korea (ando, hbkang@inha.ac.kr).} \and  Hyeonbae Kang\footnotemark[2]}

\maketitle

\begin{abstract}
We investigate in a quantitative way the plasmon resonance at eigenvalues and the essential spectrum (the accumulation point of eigenvalues) of the Neumann-Poincar\'e operator on smooth domains. We first extend the symmetrization principle so that the single layer potential becomes a unitary operator from $H^{-1/2}$ onto $H^{1/2}$. We then show that the resonance at the essential spectrum is weaker than that at eigenvalues. It is shown that anomalous localized resonance occurs at the essential spectrum on ellipses, but cloaking does not occur on ellipses unlike the core-shell structure considered in \cite{MN_06}. It is shown that resonance does not occur at the essential spectrum on three dimensional balls.
\end{abstract}



\section{Introduction}\label{sec:intro}

The Neumann-Poincar\'e (NP) operator is a boundary integral operator which arises naturally when solving classical Dirichlet and Neumann boundary value problems using layer potentials. This operator can be realized as a self-adjoint operator using Plemelj's symmetrization principle (see the next section). If the boundary of the domain is smooth, the operator is compact (see \cite{Fo95, Ke29}) and its spectrum consists of the point spectrum (eigenvalues) and the essential spectrum which is an accumulation point of the eigenvalues. The purpose of this paper is to investigate resonance at eigenvalues and at the essential spectrum, and compare them in a quantitatively precise way. The resonance at the eigenvalues of the NP operator is the plasmon resonance \cite{Grieser}. We show that the resonance at the essential spectrum (on ellipses) is the anomalous localized resonance which was first discovered on a concentric core-shell structure in \cite{NMM_94}.

To be more precise, suppose that a bounded simply connected domain $\GO$ in $\Rbb^d$ ($d=2,3$) is occupied with a plasmonic material of negative dielectric constant. In general the material property of the domain $\GO$ is represented by $\Ge_c+i \Gd$ where $\Ge_c <0$ is the dielectric constant and $\Gd >0$ indicates the dissipation. Let $\Ge_m>0$ be the dielectric constant of the matrix $\Rbb^d \setminus \overline{\GO}$. So the distribution of the dielectric constant of the structure is given by
\beq
\Ge=
\begin{cases}
\Ge_c+i\Gd \quad&\mbox{in } \GO, \\
\Ge_m \quad&\mbox{in } \Rbb^d \setminus \overline{\GO}.
\end{cases}
\eeq
We quantify the resonance through the following equation:
\beq\label{trans}
\left\{
\begin{array}{ll}
\nabla \cdot \Ge \nabla u = a \cdot \nabla \Gd_z \quad &\mbox{in } \Rbb^d, \\
u(x) = O(|x|^{1-d}) \quad &\mbox{as } |x| \to \infty,
\end{array}
\right.
\eeq
where $a$ is a constant vector and $\Gd_z$ is the Dirac mass at $z \in \Rbb^d \setminus \overline{\GO}$.
If $u_\Gd$ is the solution to \eqnref{trans}, the resonance is characterized by the blow-up of $\| \nabla u_\Gd \|_{L^2(\GO)}$:
\beq
\| \nabla u_\Gd \|_{L^2(\GO)} \to \infty \quad\mbox{as } \Gd \to 0.
\eeq
We are particularly interested in the blow-up rate of $\| \nabla u_\Gd \|_{L^2(\GO)}$ in terms of $\Gd$ when the resonance takes place.
We emphasize that as $\Gd \to 0$ $\Ge$ tends to $\Ge_c <0$ in $\GO$, and so the problem is not an elliptic one. Resonance is connected to the spectrum
of the NP operator because of the transmission conditions (continuity of the potential and the flux) along the interface $\p\GO$:
\beq\label{transcon}
u_\Gd |_+=u_\Gd|_-, \quad \Ge_m \pd{u_\Gd}{\nu} \Big|_+ = (\Ge_c+i\Gd) \pd{u_\Gd}{\nu} \Big|_- \quad\mbox{on } \p\GO,
\eeq
where the subscripts $+$ and $-$ denote limits to $\p\GO$ from outside and inside $\GO$, respectively. We will make this connection clear in section \ref{sec3}.

The findings of this paper show that the generic rate of the resonance at the eigenvalues is $\Gd^{-1}$ while that at the essential spectrum is weaker than $\Gd^{-1}$. We then show that resonance occurs at the essential spectrum ($0$) on ellipses and exact rate of resonance is provided. It turns out that it is anomalous localized resonance (see subsection \ref{sec:ell} for precise statements). We also show that resonance does not occur at the essential spectrum on the three dimensional balls. For the purpose of analysis on resonance we extend the symmetrization principle and show as a result that the single layer potential is a unitary operator from $H^{-1/2}(\p{\GO})$ onto $H^{1/2}(\p{\GO})$ ($H^s$ is a Sobolev space). We also derive an expansion formula for the fundamental solution of the Laplacian in terms of eigenfunctions of the NP operator which seems of independent interest. It is worth mentioning that it is recently found in \cite{NN_14} that on a disk a complete resonance occurs at $0$. This happens because $0$ is the only eigenvalue of the NP operator on a disk.

Recently there is rapidly growing interest in the spectral theory of the NP operator in relation to plasmonics. In \cite{KhPuSh07} the Poincar\'e's variational program was revisited with modern prospective and symmetrization of the NP operator was proved.  Eigenvalues on disks, ellipses, and balls were computed \cite{AKL, KS99}. Asymptotics of eigenvalues of the NP operator associated with closely located two-dimensional convex domain has been obtained in \cite{BT, BT2}. There has been progress on spectral theory of the NP operator on non-smooth domains. A bound on the essential spectrum on the two-dimensional curvilinear domains has been obtained \cite{PP}. Quite recently, the complete spectral resolution of the NP operator on intersecting disks has been derived which in particular shows that only the absolutely continuous spectrum exists \cite{KLS}. Interestingly, the results in \cite{KLS} shows that the spectral bound obtained in \cite{PP} is optimal in the case of intersecting disks. The spectral theory of the NP operator has been applied to analysis of cloaking by anomalous localized resonance on the plasmonic structure \cite{ACKLM} and to the study of uniformity of elliptic estimates \cite{KKLS}.

This paper is organized as follows. In section \ref{sec2} we review and extend the symmetrization of the NP operator. In section \ref{sec3} we derive an expansion formula for the fundamental solution of the Laplacian in terms of the eigenfunctions of the NP operator and obtain a representation of the solution to \eqnref{trans}. Resonance at eigenvalues and at essential spectrum is studied in section \ref{sec4} and section \ref{sec5}, respectively.
In Appendix we derive explicit expansion formula for the fundamental solution in terms of spherical harmonics and elliptic harmonics.

\section{Neumann-Poincar\'e operator and symmetrization}\label{sec2}

Let us first fix some notation.
\begin{itemize}
\item We denote by $H^{-1/2}(\p{\GO})$ the dual space of $H^{1/2}(\p{\GO})$, and by $\la \cdot , \cdot \ra$ the duality  pairing of $H^{-1/2}$ and $H^{1/2}$, and $\| \cdot \|_s$ denotes the $H^s$ norm on $\p\GO$. Let $H_0^{-1/2}(\p{\GO})$ be the space of $\psi \in H^{-1/2}(\p{\GO})$ satisfying $\la \psi, 1 \ra =0$.
\item The notation $A \lesssim B$ means that $A \le CB$ for some constant $C$, and $A \approx B$ means that both $A \lesssim B$ and $B \lesssim A$ hold.
\item Let $f (\Gd)$ and $g(\Gd)$ be positive quantities depending on $\Gd$. We write $f (\Gd) \sim g (\Gd)$ as $\Gd \to 0$ if there are constants $C_1$ and $C_2$ such that
\beq
C_1 < \frac{f (\Gd)}{g (\Gd)} < C_2 .
\eeq
\end{itemize}

Let $\GG(x)$ be the fundamental solution to the Laplacian on $\mathbb{R}^d$ ($d=2,3$), {\it i.e.},
\[
\GG(x) =
\begin{cases}
  \ds \frac{1}{2 \pi} \ln |x|, & d = 2, \\[0.2cm]
  \ds -\frac{1}{4\pi} |x|^{-1}, & d = 3.
\end{cases}
\]
The Neumann-Poincar\'e (NP) operator on $\p\GO$, denoted by $\Kcal_{\p\GO}$, is defined by
\beq\label{NPoperator}
\Kcal_{\p\GO} [\Gvf] (x) := \int_{\p \GO} \frac{\p}{\p\nu_x} \GG (x-y) \Gvf (y) \, d\Gs(y) \;, \quad x \in \p\GO.
\eeq
where $\frac{\p}{\p\nu_x}$ indicates the outward normal derivative in $x$-variable. The adjoint operator $\Kcal_{\p\GO}^*$ on $L^2(\p\GO)$ will be called the adjoint NP operator.

Importance of the adjoint operator in dealing with interface problems lies in its relation with the single layer potential $\Scal_{\p\GO}[\Gvf]$ defined by
\beq
\Scal_{\p\GO} [\Gvf] (x) := \int_{\p \GO} \Gamma (x-y) \Gvf (y) \, d\sigma(y) \;, \quad x \in \Rbb^d .
\eeq
The following jump relation holds:
\beq\label{singlejump}
\frac{\p}{\p\nu} \Scal_{\p\GO} [\Gvf] \big |_\pm (x) = \Bigl( \pm \frac{1}{2} I + \Kcal_{\p\GO}^* \Bigr) [\Gvf] (x),
\quad x \in \p\GO\;,
\eeq
where the subscripts $\pm$ indicate the limits (to $\p\GO$) from outside and inside of $\GO$, respectively.

The NP operator is not self-adjoint on the usual $L^2$-space, unless the domain $\GO$ is a disc or a ball (see \cite{Li01}). However, it is found in \cite{KhPuSh07} (see also \cite{Kang}) that the adjoint NP operator $\Kcal_{\p\GO}^*$ can be symmetrized using Plemelj's symmetrization principle (also known as Calder\'{o}n's identity):
\beq
\Scal_{\p\GO} \Kcal_{\p\GO}^* = \Kcal_{\p\GO} \Scal_{\p\GO}. \label{eq3:2}
\eeq
Define, for $\Gvf, \psi \in H^{-1/2}(\p{\GO})$,
\beq
\la \Gvf, \psi \ra_{\Hcal^*} := - \la \Gvf, \Scal_{\p\GO}[\psi] \ra. \label{eq3:3}
\eeq
Since $\Scal_{\p\GO}$ maps $H^{-1/2}(\p{\GO})$ into $H^{1/2}(\p{\GO})$, the right hand side of \eqref{eq3:3} is well-defined.
It is known that $\la \cdot, \cdot\ra_{\Hcal^*}$ is an inner product on $H_0^{-1/2}(\p{\GO})$, and the norm $\| \cdot \|_{\Hcal^*}$ induced by this inner product is equivalent to the $H^{-1/2}(\p\GO)$ norm, namely,
\beq
\| \Gvf \|_{\Hcal^*} \approx \| \Gvf \|_{-1/2}
\eeq
for all $\Gvf \in H_0^{-1/2}(\p{\GO})$ (see \cite{KKLS}).
Let $\Hcal^*_0$ be the space $H_0^{-1/2}(\p{\GO})$ equipped with the inner product $\la \cdot, \cdot\ra_{\Hcal^*}$. Then the symmetrization principle \eqref{eq3:2} shows that $\Kcal_{\p\GO}^*$ is self-adjoint on $\Hcal^*_0$.

Let us now consider symmetrization of $\Kcal_{\p\GO}$. The NP operator $\Kcal_{\p\GO}$ can be symmetrized using \eqnref{eq3:2} expressed in a different form:
\beq
\Kcal_{\p\GO}^* \Scal_{\p\GO}^{-1} = \Scal_{\p\GO}^{-1} \Kcal_{\p\GO}, \label{eq3:22}
\eeq
provided that $\Scal_{\p\GO}^{-1}$ exists. However, it was proved in \cite{Ve84} that $\Scal_{\p\GO}^{-1}$ exists only in three dimensions (or higher), and there are domains $\GO$ in two dimensions where $\Scal_{\p\GO}^{-1}$ does not exist. Here we present a simple way to overcome this difficulty.

To symmetrize $\Kcal_{\p\GO}$,
we use, as a replacement of $\Scal_{\p\GO}$ in \eqnref{eq3:2}, the operator $\Acal: H^{-1/2}(\p{\GO}) \times \mathbb{C} \to H^{1/2}(\p{\GO}) \times \mathbb{C}$ defined by
\begin{equation}\label{eq3:1}
  \Acal(\psi, a) := (\Scal_{\p\GO}[\psi] + a, \la \psi, 1 \ra).
\end{equation}
It is proved in \cite[Theorem 2.13]{AmKa07Book2} and \cite{KKLS} that $\Acal$ is invertible for $d \ge 2$. For $f \in H^{1/2}(\p{\GO})$ let
\beq
(\psi_f, a_f)= \Acal^{-1}(f,0).
\eeq
Then $\psi_f \in H_0^{-1/2}(\p{\GO})$ and it holds that
\beq\label{Saf}
\Scal_{\p\GO}[\psi_f] + a_f =f
\eeq
Moreover, we have
\beq\label{eq3:18}
\Vert{\psi_f}\Vert_{-1/2} + \vert{a_f}\vert \approx \Vert{f}\Vert_{1/2}.
\eeq
This shows in particular that the mapping $f \mapsto a$ is a bounded linear functional on $H^{1/2}(\p\GO)$. So, there is a unique $\Gvf_0 \in H^{-1/2}(\p\GO)$ such that
\beq\label{phizero}
a_f= \la \Gvf_0, f \ra
\eeq
for all $f \in H^{1/2}(\p\GO)$. Since such a $(\psi_f, a_f)$ is unique, we see that $a_f=1$ if $f=1$, in other words,
\beq\label{one}
\la \Gvf_0, 1 \ra =1.
\eeq
Since $\Kcal_{\p\GO}[1] = 1/2$, by applying $\Kcal_{\p\GO}$ to both sides of \eqnref{Saf} we see that
$$
\Kcal_{\p\GO}\Scal_{\p\GO}[\psi_f] + \frac{1}{2}a_f = \Kcal_{\p\GO}[f].
$$
It then follows from \eqnref{eq3:2} that
$$
\Scal_{\p\GO}\Kcal_{\p\GO}^*[\psi_f] + \frac{1}{2}a_f = \Kcal_{\p\GO}[f],
$$
which implies that $a_{\Kcal_{\p\GO}[f]}= \frac{1}{2}a_f$, in other words,
$$
\la \Gvf_0, \Kcal_{\p\GO}[f] \ra = \la \Gvf_0, f \ra/2.
$$
So, we have
\beq\label{eigen}
\Kcal_{\p\GO}^*[\Gvf_0]= \frac{1}{2} \Gvf_0,
\eeq
namely, $\Gvf_0$ is an eigenfunction of $\Kcal_{\p\GO}^*$ (on $H^{-1/2}(\p\GO)$ corresponding to the eigenvalue $1/2$ normalized by \eqnref{one}. We can infer from the jump formula \eqnref{singlejump} and \eqnref{eigen} that $\Scal_{\p\GO}[\Gvf_0]$ is constant in $\GO$.
We emphasize that the function $\Gvf_0$ already appeared in literatures. It is proved in \cite{Ve84} that $\Scal_{\p\GO}[\Gvf_0]=0$ in $\GO$ for some domain $\GO$ in two dimensions, which is why $\Scal_{\p\GO}$ is not invertible.

We now define a variant of the single layer potential on $\Hcal^*$ by
\beq
\widetilde{\Scal}_{\p\GO}[\Gvf]=
\begin{cases}
\Scal_{\p\GO}[\Gvf] \quad &\mbox{if } \la \Gvf, 1 \ra=0, \\
1 \quad &\mbox{if } \Gvf=\Gvf_0.
\end{cases}
\eeq
Then we have an extension of \eqnref{eq3:2}:
\beq
\widetilde{\Scal}_{\p\GO} \Kcal_{\p\GO}^* = \Kcal_{\p\GO} \widetilde{\Scal}_{\p\GO}.
\eeq
Using this we can extend the inner product \eqnref{eq3:3} (defined on $H^{-1/2}_0(\p\GO)$) to $H^{-1/2}(\p\GO)$. For $\Gvf, \psi\in H^{-1/2}(\p\GO)$, define
\beq\label{ext}
\la \Gvf, \psi \ra_{\Hcal^*} := - \la \Gvf, \widetilde{\Scal}_{\p\GO}[\psi] \ra .
\eeq
We also have a new inner product on $H^{1/2}(\p\GO)$:
\beq\label{ext2}
\la \Gvf, \psi \ra_{\Hcal} := - \la \Gvf, \widetilde{\Scal}_{\p\GO}^{-1}[\psi] \ra .
\eeq
If we define $\Hcal$ to be $H^{1/2}(\p\GO)$ with this inner product, then $\Kcal_{\p\GO}$ is self-adjoint on $\Hcal$. We emphasize that norm $\| \cdot \|_{\Hcal}$ is equivalent to $\| \cdot \|_{1/2}$. Moreover, $\widetilde{\Scal}_{\p\GO}: \Hcal^* \to \Hcal$ is a unitary operator. Observe that if $\{ \psi_j\}_{j=1}^\infty$ is an orthonormal basis of $\Hcal^*_0$, then $\{ \psi_j\}_{j=1}^\infty \cup \{ \Gvf_0 \}$ is an orthonormal basis of $\Hcal^*$ and $\{ \Scal_{\p\GO}[\psi_j]\}_{j=1}^\infty \cup \{ 1 \}$ is an orthonormal basis of $\Hcal$.

\section{Representation of the solution}\label{sec3}

If $\p \GO$ is $\Ccal^{1, \Ga}$-smooth for some $\Ga>0$, then it is known that $\Kcal_{\p\GO}^*$ is a compact operator on $\Hcal^*$ (see \cite{Ke29}). Since $\Kcal_{\p\GO}^*$ is self-adjoint on $\Hcal^*_0$, its eigenvalues $\{\Gl_j\}_{j=1}^{\infty}$ are real and accumulates to $0$. We emphasize that $|\Gl_j| < 1/2$ (see \cite{ChLe08, Ve84}).
Let $\{\psi_j\}_{j=1}^{\infty}$ with $\| \psi_j \|_{\Hcal^*}=1$ be the corresponding (real valued) eigenfunctions counting the multiplicities. Then we have shown in the previous section that $\{ \Scal_{\p\GO}[\psi_j]\}_{j=1}^\infty \cup \{ 1 \}$ is an orthonormal basis of $\Hcal$.

Fix $z \in \Rbb^d \setminus\overline{\GO}$. Then $\GG(\cdot-z)$ belongs to $H^{1/2}(\p\GO)$, and so admits the following decomposition:
\beq
\GG(x-z) = \sum_{j=1}^\infty c_j(z) \Scal_{\p\GO}[\psi_j](x) + c_0(z), \quad x \in \p\GO,
\eeq
for some constants $c_j(z)$ (depending on $z$) satisfying
\beq\label{cjfinite}
\sum_{j=1}^\infty |c_j(z)|^2 < \infty.
\eeq
Since $- \la \Scal_{\p\GO}[\psi_j], \psi_i \ra= \Gd_{ij}$ where $\Gd_{ij}$ is the Kronecker's delta, we see that
$$
c_j(z) =- \int_{\p\GO} \GG(x-z) \psi_j(x) d\Gs(x)= - \Scal_{\p\GO}[\psi_j](z), \quad j=1,2,3,\ldots.
$$
We also see from \eqnref{one} that
$$
c_0(z) = \Scal_{\p\GO}[\Gvf_0](z), \quad j=1,2,3,\ldots.
$$
So, we obtain the following formula:
\beq
\GG(x-z) = - \sum_{j=1}^\infty \Scal_{\p\GO}[\psi_j](z) \Scal_{\p\GO}[\psi_j](x) + \Scal_{\p\GO}[\Gvf_0](z), \quad x \in \p\GO,
\eeq
Observe that
\beq\label{localized}
\left\| \sum_{j=1}^\infty \Scal_{\p\GO}[\psi_j](z) \Scal_{\p\GO}[\psi_j] \right\|_{\Hcal}^2 =  \sum_{j=1}^\infty |\Scal_{\p\GO}[\psi_j](z)|^2 < \infty.
\eeq
Since $\| \cdot \|_{\Hcal}$ is equivalent to $\| \cdot \|_{1/2}$ as shown in the previous section, we find from the trace theorem that $\sum_{j=1}^\infty \Scal_{\p\GO}[\psi_j](z) \Scal_{\p\GO}[\psi_j]$ converges in $H^1(\GO)$ and is harmonic in $\GO$. So, we obtain the following theorem
on an expansion of the fundamental solution to the Laplacian.

\begin{theorem}\label{thm:fund}
It holds that
\beq\label{decomp}
\GG(x-z) = - \sum_{j=1}^\infty \Scal_{\p\GO}[\psi_j](z) \Scal_{\p\GO}[\psi_j](x) + \Scal_{\p\GO}[\Gvf_0](z), \quad x \in \overline{\GO}, \ z \in \Rbb^d \setminus\overline{\GO}.
\eeq
\end{theorem}

If $\GO$ is a ball, then $\Scal_{\p\GO} [\psi_j]$ is a spherical harmonics (see subsection \ref{sec:ball}). Therefore \eqnref{decomp} is the expansion of $\GG(x-z)$ in terms of the spherical harmonics, which is well-known (see, for example, \cite{Ke29}). We will show explicit formula for \eqnref{decomp} when $\GO$ is a unit ball in $\Rbb^3$ or an ellipse in $\Rbb^2$ in the Appendix.
The formula \eqnref{decomp} shows that the fundamental solution can be expanded on any smooth domain in terms of the localized plasmon (the single layer potential of an eigenfunction is called a localized plasmon). The estimate \eqnref{localized} shows that $\Scal_{\p\GO}[\psi_j](z) \to 0$ as $j \to \infty$ for all $z \notin \overline{\GO}$. So, if $j$ is large then $\Scal_{\p\GO}[\psi_j]$ is localized near $\GO$. It explains why $\Scal_{\p\GO}[\psi_j]$ is called a {\it localized} plasmon.

We now derive a representation of the solution to \eqnref{trans}. Let
\beq\label{Fzx}
F_z(x) := -a \cdot \nabla_x \GG (x-z), \quad x \neq z.
\eeq
Then $\GD F_z(x) = a \cdot \nabla \Gd_z(x)$, and hence we see from transmission conditions \eqnref{transcon} on $\p\GO$ and \eqnref{singlejump} that the solution $u_\Gd$ to \eqnref{trans} can be represented as
\beq\label{solrep}
u_\Gd(x) = F_z(x) + \Scal_{\p\GO} [\Gvf_\Gd](x), \quad x \in \Rbb^d,
\eeq
where the potential $\Gvf_\Gd \in \Hcal^*_0(\p\GO)$ is the solution to the integral equation
\beq\label{solint}
\left( \Gl I - \Kcal_{\p\GO}^* \right)[\Gvf_\Gd] = \p_\nu F_z \quad\mbox{on } \p\GO
\eeq
($\p_\nu F_z$ denotes the normal derivative of $F_z$). Here
\beq\label{lambda}
\Gl:= \frac{\Ge_c +\Ge_m + i\Gd}{2(\Ge_c-\Ge_m) + 2i\Gd}.
\eeq
Note that $\Gl \to \frac{\Ge_c +\Ge_m}{2(\Ge_c-\Ge_m)}$ as $\Gd \to 0$. The number $\Ge_c/\Ge_m$ such that
$\frac{\Ge_c+\Ge_m}{2(\Ge_c-\Ge_m)}$ is an eigenvalue of $\Kcal_{\p\GO}^*$ is called a plasmonic eigenvalue \cite{Grieser}. We emphasize that $\Ge_c$ is negative if and only if  $\frac{\Ge_c+\Ge_m}{2(\Ge_c-\Ge_m)}$ lies in $(-1/2, 1/2)$ where the spectrum of $\Kcal_{\p\GO}^*$ lies.

Since $\Kcal_{\p\GO}^*$ admits the spectral decomposition
\beq\label{Expansion_K}
\Kcal_{\p\GO}^* = \sum_{j = 1}^{\infty} \Gl_j \psi_j \otimes \psi_j,
\eeq
the solution $\Gvf_\Gd$ to \eqnref{solint} can be expressed as
\beq\label{Expansion_phi_j}
\Gvf_\Gd = \sum_{j = 1}^{\infty} \frac{\Ga_j(z)}{\Gl - \Gl_j}  \psi_j,
\eeq
where
\beq
\Ga_j(z):= \left\la \p_\nu F_z, \psi_j \right\ra_{\Hcal^*}.
\eeq
We can see from \eqnref{Fzx} that
$$
\Ga_j(z) = -a \cdot \nabla \int_{\p\GO} \pd{}{\nu_x} \GG(x-z) \Scal_{\p\GO}[\psi_j](x) \, d\Gs(x) .
$$
It can seen from \eqnref{singlejump} and \eqnref{decomp} that
$$
\pd{}{\nu_x} \GG(x-z) = - \sum_{j=1}^\infty \Scal_{\p\GO}[\psi_j](z) \pd{}{\nu} \Scal_{\p\GO}[\psi_j](x) =
\sum_{j=1}^\infty \Big( \frac{1}{2} - \Gl_j \Big) \Scal_{\p\GO}[\psi_j](z) \psi_j(x).
$$
It then follows that
\beq\label{fj}
\Ga_j(z) = \Big( \Gl_j - \frac{1}{2} \Big) a \cdot \nabla \Scal_{\p\GO}[\psi_j](z).
\eeq

\section{Resonance at eigenvalues}\label{sec4}

We investigate the behavior of the solution $u_\Gd$ when $\Gl$ approaches to one of eigenvalues $\Gl_l \neq 0$ as $\Gd \to 0$, namely, when
\beq
\frac{\Ge_c +\Ge_m}{2(\Ge_c-\Ge_m)}= \Gl_l.
\eeq
We show that
\beq\label{estimate41}
\| \nabla u_\Gd \|_{L^2(\GO)} \sim \Gd^{-1} \quad\mbox{as } \Gd \to 0,
\eeq
as one may expect.

We first show that
\beq\label{gradest}
\| \nabla {\Scal}_{\p\GO}[\Gvf] \|_{L^2(\GO)}^2 \approx \| \Gvf \|_{\Hcal^*}^2
\eeq
for all $\Gvf \in \Hcal^*_0$. In fact, we see from \eqnref{singlejump} and \eqnref{Expansion_K} that
\begin{align*}
\| \nabla {\Scal}_{\p\GO}[\Gvf] \|_{L^2(\GO)}^2 &= \int_{\p\GO} {\Scal}_{\p\GO}[\Gvf] \overline{\pd{}{\nu} {\Scal}_{\p\GO}[\Gvf] \Big|_{-}} \, d\Gs \\
&= \left\la \Gvf, (-\frac{1}{2}I + \Kcal_{\p\GO}^*)[\Gvf] \right\ra_{\Hcal^*} \\
&= \sum_{j=1}^\infty (\frac{1}{2}- \Gl_j) |\la \Gvf, \psi_j \ra_{\Hcal^*}|^2.
\end{align*}
Since $|\Gl_j| < 1/2$ and accumulates to $0$, we have \eqnref{gradest}.

We now see from \eqnref{Expansion_phi_j} that
$$
\| \nabla (u_\Gd - F_z) \|_{L^2(\GO)}^2 \approx \| \Gvf_\Gd \|_{\Hcal^*}^2 =  \sum_{\Gl_j=\Gl_l} \frac{|\Ga_j(z)|^2}{|\Gl - \Gl_l|^2} + \sum_{\Gl_j \neq \Gl_l} \frac{|\Ga_j(z)|^2}{|\Gl - \Gl_j|^2}.
$$
If $\Gl_j \neq \Gl_l$, then $|\Gl_j - \Gl_l| > C$ for some positive constant $C$ since $\Gl_l \neq 0$, and so the second term on the righthand side above is bounded as $\Gd \to 0$. Since $|\Gl-\Gl_l| \sim \Gd$ as $\Gd \to 0$, we obtain
$$
\| \nabla (u_\Gd - F_z) \|_{L^2(\GO)} \sim \Gd^{-1}.
$$
Since $\| \nabla F_z \|_{L^2(\GO)}$ is bounded, \eqnref{estimate41} follows.

\section{Resonance at the essential spectrum}\label{sec5}

In this section we assume that $0$ is not an eigenvalue of $\Kcal_{\p\GO}^*$. Since eigenvalues of $\Kcal_{\p\GO}^*$ converges to $0$, $\{ 0 \}$ is the essential spectrum of $\Kcal_{\p\GO}^*$.
It is worth mentioning that we are not aware of any domain other than disks on which the NP operator has $0$ as an eigenvalue. If $\GO$ is a disk, then $\Kcal_{\p\GO}^* \equiv 0$ on $\Hcal^*_0$.

We consider the resonance when $\Gl \to 0$ as $\Gd \to 0$, in other words, when
\beq
\Ge_c + \Ge_m=0.
\eeq
In this case, we assume that $\Gl=i\Gd$ for simplicity of presentation.

We first obtain the following theorem which shows that the resonance, if it occurs, is never at the rate of $\Gd^{-1}$ unlike the resonance at eigenvalues.
\begin{theorem}
It holds that
\beq
\lim_{\Gd \to 0} \Gd \| \nabla u_\Gd \|_{L^2(\GO)} =0.
\eeq
\end{theorem}

\pf
We obtain from \eqnref{Expansion_phi_j} and \eqnref{gradest} that
\beq\label{anomal1}
\| \nabla (u_\Gd -F_z) \|_{L^2(\GO)}^2 \sim \sum_{j=1}^\infty \frac{|\Ga_j(z)|^2}{\Gd^2 + \Gl_j^2} .
\eeq
We then decompose the summation into two parts as
\beq
\sum_{j=1}^\infty \frac{|\Ga_j(z)|^2}{\Gd^2 + \Gl_j^2} = \sum_{|\Gl_j| \le \Gd} \frac{|\Ga_j(z)|^2}{\Gd^2 + \Gl_j^2} + \sum_{|\Gl_j| > \Gd} \frac{|\Ga_j(z)|^2}{\Gd^2 + \Gl_j^2} =: S_1 + S_2.
\eeq
Since $\sum |\Ga_j(z)|^2 <\infty$, we have
$$
\Gd^2 S_1 \le \sum_{|\Gl_j| \le \Gd} |\Ga_j(z)|^2 \to 0 \ \mbox{as } \Gd \to 0.
$$

To show that $\Gd^2 S_2 \to 0$, we express $S_2$ as
$$
S_2 = \sum_{k=0}^\infty \sum_{2^k \Gd < |\Gl_j| \le 2^{k+1} \Gd} \frac{|\Ga_j(z)|^2}{\Gd^2 + \Gl_j^2}.
$$
Then we see that
$$
\Gd^2 S_2 \le \sum_{k=0}^\infty \frac{1}{1 + 2^{2k}}\sum_{2^k \Gd < |\Gl_j| \le 2^{k+1} \Gd} |\Ga_j(z)|^2
\le \sum_{k=0}^\infty \frac{1}{1 + 2^{2k}}\sum_{|\Gl_j| \le 2^{k+1} \Gd} |\Ga_j(z)|^2 .
$$
For each fixed $k$, $\sum_{|\Gl_j| \le 2^{k+1} \Gd} \Ga_j(z)^2 \to 0$ as $\Gd \to 0$. So we infer that $S_2 \to 0$ as $\Gd \to 0$ by Lebesgue dominated convergence theorem.
The proof is complete. \qed

In order to derive estimates for the resonance, we need to investigate the asymptotic behavior of the quantity on the righthand side of \eqnref{anomal1} as $\Gd \to 0$. This seems a quite difficult task since it depends on the behavior of $\Gl_j$ and $\Ga_j(z)$ as $j \to \infty$. So,
we deal with two specific domains, ellipses and three dimensional balls, where eigenvalues and eigenfunctions are known. It is worth mentioning that the disks are out of consideration since the NP operator on disks are $0$.

\subsection{Anomalous localized resonance on ellipses}\label{sec:ell}

We first consider the resonance when $\GO$ is an ellipse in $\Rbb^2$. The elliptic coordinates $x = (x_1, x_2) = (x_1(\rho, \Go), x_2(\rho, \Go))$ is given by
\[
x_1 (\rho, \Go) = R \cos{\Go} \cosh{\rho}, \quad x_2 (\rho, \Go) = R \sin{\Go} \sinh{\rho}, \quad \rho > 0, \ 0 \le \Go < 2 \pi.
\]
When the above holds, we denote $\rho=\rho_x$ and $\Go_x$. Then we can represent $\p\GO$ as
\beq\label{ellipse}
\p\GO = \{ x \in \Rbb^2; \rho_x=\rho_0 \}
\eeq
for some $\rho_0 > 0$. The number $\rho_0$ is called the elliptic radius of $\GO$. The position $z$ of the source is denote by $\rho_z$ and $\Go_z$ in elliptic coordinates. We obtain the following theorem whose proof will be given in the next section.

\begin{theorem}\label{anomalellipse}
Suppose that $\GO$ is an ellipse given by \eqnref{ellipse}. Then we have
\beq\label{anomalellipse2}
  \| \nabla u_\Gd \|_{L^2(\GO)}^2 \sim
  \begin{cases}
    \Gd^{-3 + \rho_z/\rho_0} |\log {\Gd}| & \ \text{if } \rho_0 < \rho_z < 3 \rho_0, \\
    |\log {\Gd}|^2 & \ \text{if } \rho_z = 3 \rho_0, \\
    1 & \ \text{if } \rho_z > 3 \rho_0,
  \end{cases}
\eeq
as $\Gd \to 0$.
\end{theorem}

The quantity $E_\Gd:=\Gd \| \nabla u_\Gd \|_{L^2(\GO)}^2$ is of particular interest since it represents the imaginary part of the energy, namely,
$$
\Gd \| \nabla u_\Gd \|_{L^2(\GO)}^2 = \Im \int_{\Rbb^d} \Ge |\nabla u_\Gd|^2 dx.
$$
Physically it represents the electro-magnetic energy dissipated into heat. Estimates \eqnref{anomalellipse2} shows that $E_\Gd \to \infty$ if $\rho_0 < \rho_z \le 2 \rho_0$ while it tends to $0$ as $\Gd \to 0$ if $\rho_z > 2 \rho_0$. So the critical (elliptic) radius is $2\rho_0$. This phenomenon is reminiscent of the anomalous localized resonance (ALR) discovered in \cite{NMM_94, MN_06} (see also \cite{ACKLM, CKKL, KLSW}). There it is shown that a disk (a core) is coated by a concentric disk (a shell) of plasmonic material, then ALR occurs if the source is located within a critical radius, and does not occur for sources outside the radius. This is exactly what \eqnref{anomalellipse2} shows. It also in accordance with the result in \cite{CKKL} where the core-shell structure of confocal ellipses of radii $\rho_e > \rho_i$ is considered. There it is shown that
the critical (elliptic) radius is given by
$$
\rho_* =
\begin{cases}
  (3 \rho_e - \rho_i) / 2, & \text{ if } \rho_e \le 3 \rho_i, \\
  2 (\rho_e - \rho_i), & \text{ if } \rho_e > 3 \rho_i. 
\end{cases}
$$
If the core shrinks, in other words $\rho_i$ tends to $0$, then $\rho_*$ tends to $2 \rho_e$ which is the critical radius found in this paper. 

ALR also requires the solution to be bounded outside a bounded set. This requirement is satisfied on ellipse as the following theorem shows. Theorem \ref{anomalellipse} and Theorem \ref{anomalellipse3} show that ALR may occur not only on the coated structures but also on simply connected structure. So ALR may be regarded as resonance at the accumulation points of eigenvalues. We emphasize that the NP operator on the coated disk has $0$ as its essential spectrum as proved in \cite{ACKLM}.

\begin{theorem}\label{anomalellipse3}
Let \(\GO\) be an ellipse given by \eqnref{ellipse}. It holds for all $x$ satisfying $\rho_x+\rho_z- 4 \rho_0 >0$ that
\beq\label{47}
|u_{\Gd} (x) - F_z (x)| \lesssim \sum_{n = 1}^{\infty} e^{-n (\rho_x+\rho_z- 4 \rho_0)}.
\eeq
In particular, let \(\overline{\rho} > 0\) be such that \(\overline{\rho} > 4 \rho_0 - \rho_z\), then there exists some \(C = C_{\overline{\rho}} > 0\) such that
  \begin{equation}
    \label{eq:56}
    \sup_{\rho_x \ge \overline{\rho}} \left\vert u_{\Gd} (x) - F_z (x) \right\vert < C .
  \end{equation}
\end{theorem}

To prove Theorem \ref{anomalellipse} and Theorem \ref{anomalellipse3} let us recall some facts on the NP operator on ellipses.
It is proved in \cite{AKL} and \cite{CKKL} that eigenvalues of $\Kcal_{\p\GO}^*$ are
\beq
\Gl_n = \frac{1}{2e^{2 n \rho_0}} , \quad n = 1, 2, \cdots,
\eeq
and corresponding eigenfunctions are
\beq
\phi_n^c (\Go) := \Xi (\rho_0, \Go)^{-1} \cos{n \Go}, \quad \phi_n^s (\Go) := \Xi (\rho_0, \Go)^{-1} \sin{n \Go}, \quad n = 1, 2, \cdots,
\eeq
where
\beq\label{Xi}
\Xi = \Xi (\rho_0, \Go) := R \sqrt{\sinh^2{\rho_0} + \sin^2{\Go}}.
\eeq
It is also proved in the same papers that
\beq\label{singlecos}
\Scal_{\p\GO} [\phi_n^c] (x) =
\begin{cases}
- \displaystyle \frac{\cosh{n \rho}}{n}  e^{-n \rho_0}\cos{n \Go}, \quad & \rho < \rho_0, \\
- \displaystyle \frac{\cosh{n \rho_0}}{n}  e^{-n \rho}\cos{n \Go}, \quad & \rho \ge \rho_0,
\end{cases}
\eeq
and
\beq\label{singlesin}
\Scal_{\p\GO} [\phi_n^s] (x) =
\begin{cases}
- \displaystyle \frac{\sinh{n \rho}}{n}  e^{-n \rho_0}\sin{n \Go}, \quad & \rho < \rho_0, \\
- \displaystyle \frac{\sinh{n \rho_0}}{n}  e^{-n \rho}\sin{n \Go}, \quad & \rho \ge \rho_0.
\end{cases}
\eeq
Here and afterwards $(\rho, \Go)$ is the elliptic coordinates of $x$. Since the length element $d \sigma$ on $\p\GO$ is given by
\beq\label{length}
d \Gs = \Xi d \Go,
\eeq
one can see that
$$
\| \phi_n^c \|_{\Hcal^*}^2 = -\la \phi_n^c, \Scal_{\p\GO}[\phi_n^c] \ra = \frac{\pi \cosh{n \rho_0}}{n e^{n \rho_0}},
$$
and
$$
\| \phi_n^s \|_{\Hcal^*}^2 = \frac{\pi\sinh{n \rho_0}}{n e^{n \rho_0}}.
$$
So, the normalized eigenfunctions are
\beq\label{515}
\psi_n^c := \sqrt{\frac{n e^{n \rho_0}}{\pi \cosh{n \rho_0}}} \phi_n^c, \quad \psi_n^s := \sqrt{\frac{n e^{n \rho_0}}{\pi \sinh{n \rho_0}}} \phi_n^s.
\eeq
We see from \eqnref{singlecos}
\beq\label{415}
\Scal_{\p\GO} [\psi_n^c] (x) =
\begin{cases}
- \sqrt{\displaystyle \frac{e^{- n \rho_0}}{n\pi\cosh{n\rho_0}}} \cosh{n \rho}\cos{n \Go}, \quad & \rho < \rho_0, \\
- \sqrt{\displaystyle \frac{e^{n \rho_0} \cosh{n \rho_0}}{n\pi}} e^{-n \rho}\cos{n \Go}, \quad & \rho \ge \rho_0,
\end{cases}
\eeq
and from \eqnref{singlesin}
\beq\label{416}
\Scal_{\p\GO} [\psi_n^s] (x) =
\begin{cases}
- \sqrt{\displaystyle \frac{e^{- n \rho_0}}{n\pi\sinh{n\rho_0}}} \sinh{n \rho}\sin{n \Go}, \quad & \rho < \rho_0, \\
- \sqrt{\displaystyle \frac{e^{n \rho_0} \sinh{n \rho_0}}{n\pi}} e^{-n \rho}\sin{n \Go}, \quad & \rho \ge \rho_0.
\end{cases}
\eeq

\medskip
\noindent{\sl Proof of Theorem \ref{anomalellipse}}.
We see from \eqnref{fj} and \eqnref{anomal1} that
\begin{align*}
& \| \nabla (u_\Gd -F_z) \|_{L^2(\GO)}^2 \sim \sum_{j=1}^\infty \frac{|\Ga_j(z)|^2}{\Gd^2 + \Gl_j^2} \\
    =  & \sum_{n = 1}^{\infty} \frac{1}{\Gd^2 + \Gl_n^2} \left[ \left\vert a \cdot \nabla \Scal_{\p\GO} [\psi_n^c] (z) \right\vert^2 + \left\vert a \cdot \nabla \Scal_{\p\GO} [\psi_n^s] (z) \right\vert^2 \right]  \\
    = & \sum_{n = 1}^{\infty} \frac{1}{\Gd^2 + \Gl_n^2} \Big[ \frac{e^{n \rho_0} \cosh{n \rho_0}}{n \pi} \left\vert a \cdot \nabla_z \left( e^{- n \rho_z} \cos {n \Go_z} \right) \right\vert^2 + \frac{e^{n \rho_0} \sinh {n \rho_0}}{n \pi} \left\vert a \cdot \nabla_z \left( e^{- n \rho_z} \sin {n \Go_z}\right) \right\vert^2 \Big].
\end{align*}
Since $\cosh{n \rho_0} \approx \sinh {n \rho_0} \approx e^{n \rho_0}$, we see that
\beq\label{57}
\sum_{j = 1}^{\infty} \frac{|\Ga_j(z)|^2}{\Gd^2 + \Gl_j^2} \sim
    \sum_{n = 1}^{\infty} \frac{1}{\Gd^2 + \Gl_n^2} \cdot \frac{e^{2n \rho_0}}{n} \Big[ \left\vert a \cdot \nabla_z \left( e^{- n \rho_z} \cos {n \Go_z} \right) \right\vert^2 \\
    + \left\vert a \cdot \nabla_z \left( e^{- n \rho_z} \sin {n \Go_z}\right) \right\vert^2 \Big]
\eeq
as $\Gd \to 0$.

Let $U(\Go)$ be the rotation by the angle $\Go$, namely,
$$
U(\Go)= \begin{bmatrix}
      \cos {\Go} & - \sin {\Go} \\
      \sin {\Go} & \cos {\Go}
    \end{bmatrix}.
$$
Using the change of variables formula
$$
\frac{\p}{\p x_1} = \frac{R}{\Xi (\rho, \Go)^2}\left( \cos {\Go} \sinh {\rho} \frac{\p}{\p \rho} - \sin {\Go} \cosh {\rho} \frac{\p}{\p \Go} \right)
$$
and
$$
\frac{\p}{\p x_2} = \frac{R}{\Xi (\rho, \Go)^2}\left( \sin {\Go} \cosh {\rho} \frac{\p}{\p \rho} + \cos {\Go} \sinh {\rho} \frac{\p}{\p \Go} \right)
$$
where $\Xi (\rho, \Go)$ is given by \eqnref{Xi}, we see through tedious but straightforward computations (we omit the computations) that
$$
a \cdot \nabla \left( e^{- n \rho} \cos {n \Go} \right)
    = \frac{- R n e^{-n \rho}}{\Xi (\rho, \Go)^2} a \cdot U(n \Go) b (\rho, \Go)
$$
and
$$
a \cdot \nabla \left( e^{- n \rho} \sin {n \Go} \right)
    = \frac{- R n e^{-n \rho}}{\Xi (\rho, \Go)^2} a \cdot U(n \Go - \pi / 2) b (\rho, \Go),
$$
where
\[
b (\rho, \Go) = \begin{bmatrix} \cos {\Go} \sinh {\rho} \\ \sin {\Go} \cosh {\rho} \end{bmatrix}.
\]
Let $\theta_n$ be the angle between the vectors $a$ and $U(n \Go) b (\rho, \Go)$. Then we have
$$
|a \cdot U(n \Go) b (\rho, \Go)|^2 = |a|^2 |b (\rho, \Go)|^2 \cos^2 {\theta_n}
$$
and
$$
|a \cdot U(n \Go - \pi/2) b (\rho, \Go)|^2 = |a|^2 |b (\rho, \Go)|^2 \sin^2 {\theta_n},
$$
which implies
\beq\label{418}
\left\vert a \cdot \nabla \left( e^{- n \rho} \cos {n \Go} \right) \right\vert^2 + \left\vert a \cdot \nabla \left( e^{- n \rho} \sin {n \Go} \right) \right\vert^2 = \frac{n^2 e^{- 2 n \rho} \left\vert a \right\vert^2 \left\vert b (\rho, \Go) \right\vert^2}{R^2 \left( \sinh^2 {\rho} + \sin^2 {\Go} \right)^2}.
\eeq
It then follows from \eqnref{57} that
\beq\label{58}
\sum_{j = 1}^{\infty} \frac{|\Ga_j(z)|^2}{\Gd^2 + \Gl_j^2} \sim
\sum_{n = 1}^{\infty} \frac{n e^{2n \rho_0} e^{- 2 n \rho_z}}{\Gd^2 + \Gl_n^2}
\eeq
as $\Gd \to 0$.

We now investigate the asymptotic behavior of the righthand side of \eqnref{58} as $\Gd \to 0$.
Let
\beq
N= \left[ -\frac{1}{2\rho_0} \log 2\Gd \right],
\eeq
which is the first number such that $\Gd > \frac{1}{2} e^{-2N \rho_0}$.
Then one can easily see that
$$
\sum_{n = 1}^{\infty} \frac{n e^{2n \rho_0} e^{- 2 n \rho_z}}{\Gd^2 + \Gl_n^2} = \sum_{n \le N} + \sum_{n > N} \sim
\sum_{n \le N} \frac{n e^{2n \rho_0} e^{- 2 n \rho_z}}{e^{-4n \rho_0}}
+ \frac{1}{\Gd^2} \sum_{n > N} n e^{-2n (\rho_z-\rho_0)} .
$$
Observe that
$$
\sum_{n \le N} \frac{n e^{2n \rho_0} e^{- 2 n \rho_z}}{e^{-4n \rho_0}} \sim
\sum_{n \le N} n e^{2n (3\rho_0-\rho_z)} \sim
\begin{cases}
1 \quad&\mbox{if } \rho_z > 3 \rho_0, \\
|\log\Gd|^2 \quad&\mbox{if } \rho_z = 3 \rho_0, \\
|\log\Gd| \Gd^{-3+ \rho_z/\rho_0} \ &\mbox{if } \rho_0 < \rho_z < 3 \rho_0.
\end{cases}
$$
On the other hand, we have
$$
\frac{1}{\Gd^2} \sum_{n > N} n e^{-2n (\rho_z-\rho_0)} \sim |\log\Gd| \Gd^{-3 + \rho_z/\rho_0}.
$$
So we infer that
$$
\sum_{n = 1}^{\infty} \frac{n e^{2n \rho_0} e^{- 2 n \rho_z}}{\Gd^2 + \Gl_n^2} \sim
\begin{cases}
1 \quad&\mbox{if } \rho_z > 3 \rho_0, \\
|\log\Gd|^2 \ &\mbox{if } \rho_z = 3 \rho_0, \\
|\log\Gd| \Gd^{-3+ \rho_z/\rho_0} \ &\mbox{if } \rho_0 < \rho_z < 3 \rho_0.
\end{cases}
$$
Since $\| \nabla F_z \|_{L^2(\GO)}^2$ is bounded, we obtain \eqnref{anomalellipse2}.  \qed

\medskip
\noindent{\sl Proof of Theorem \ref{anomalellipse3}}.
One can see from \eqnref{solrep}, \eqnref{Expansion_phi_j} and \eqnref{fj} that
\begin{align*}
u_{\Gd} (x) - F_z (x) =  & \sum_{n = 1}^{\infty} \frac{1}{i\Gd - \Gl_n} \left( \frac{1}{2} - \Gl_n \right)
\Big[ \left( a \cdot \nabla_z \Scal_{\p\GO} [\psi_n^c] (z) \right) \Scal_{\p\GO} [\psi_n^c] (x)  \\
    & \quad \quad \ + \left( a \cdot \nabla_z \Scal_{\p\GO} [\psi_n^s] (z) \right) \Scal_{\p\GO} [\psi_n^s] (x) \Big].
\end{align*}
It then follows from \eqnref{415} and \eqnref{416} that
\begin{align*}
u_{\Gd} (x) - F_z (x) =
    & \sum_{n = 1}^{\infty} \frac{1}{i\Gd - \Gl_n} \left( \frac{1}{2} - \Gl_n \right)  \\
    & \quad \times \Big [\frac{e^{n \rho_0} \cosh {n \rho_0}}{n \pi} \left( a \cdot \nabla_z \left( e^{- n \rho_z} \cos {n \Go_z} \right) \right) e^{- n \rho_x} \cos {n \Go_x} \nonumber \\
    & \quad \quad + \frac{e^{n \rho_0} \sinh {n \rho_0}}{n \pi} \left( a \cdot \nabla_z \left( e^{- n \rho_z} \sin {n \Go_z} \right) \right) e^{- n \rho_x} \sin {n \Go_x} \Big]
\end{align*}
where $(\rho_z, \Go_z)$ is the elliptic coordinates of $z$. Therefore we have
\begin{align*}
|u_{\Gd} (x) - F_z (x)| \lesssim
    & \sum_{n = 1}^{\infty} \frac{e^{4n \rho_0}}{n}
    \Big [\left| a \cdot \nabla_z \left( e^{- n \rho_z} \cos {n \Go_z} \right) \right| + \left| a \cdot \nabla_z \left( e^{- n \rho_z} \sin {n \Go_z} \right) \right| \Big]  e^{- n \rho_x}
\end{align*}
We then see from \eqnref{418} that
$$
|u_{\Gd} (x) - F_z (x)| \lesssim \sum_{n = 1}^{\infty} \frac{e^{4n \rho_0}}{n}
n e^{-  n \rho_z} e^{- n \rho_x} = \sum_{n = 1}^{\infty} e^{-n (\rho_x+\rho_z- 4 \rho_0)}.
$$
\eqnref{eq:56} is an immediate consequence of \eqnref{47}. This completes the proof. \qed

\subsection{Anomalous localized resonance on balls}\label{sec:ball}

The following theorem shows that resonance at the essential spectrum does not occur if $\GO$ is a three dimensional ball.
\begin{theorem}\label{anomalball}
Suppose that $\GO$ is a three dimensional ball and $f(x)= a \cdot \nabla \Gd_z(x)$  for some $z \in \Rbb^3 \setminus \overline{\GO}$. Then there is a constant $C$ such that
\beq\label{anomalball2}
  \| \nabla u_\Gd \|_{L^2(\GO)} \le C
\eeq
for all $\Gd$.
\end{theorem}

\pf
Assume that \(\GO\) is the unit ball centered at $0$ for convenience. Then eigenvalue of $\Kcal_{\p\GO}^*$ are
\beq
\Gl_n = \frac{1}{2 (2 n + 1)}, \quad n=1,2, \ldots,
\eeq
and corresponding eigenfunctions are \(Y_n^m (\hat{x})\), $m = - n, - n + 1, \cdots, n - 1, n$, the spherical harmonics of degree \(n\), where \(\hat{x} = x / | x |\) (see \cite{ACKLM14}). It is also proven in the same paper that
\beq\label{singleball}
\Scal_{\p\GO} [Y_n^m] (x) =
\begin{cases}
- \displaystyle \frac{1}{2 n + 1} r^n Y_n^m (\hat{x}), \quad & \text{ for } |x| = r < 1, \\
- \displaystyle \frac{1}{2 n + 1} r^{-(n + 1)} Y_n^m (\hat{x}), \quad & \text{ for } |x| = r \ge 1.
\end{cases}
\eeq
So the eigenfunctions normalized in $\Hcal^*$ are
\beq\label{eigenfunball}
\psi_n^m (x) = \frac{1}{\sqrt{2 n + 1}} Y_n^m (x), \quad m = - n, - n + 1, \cdots, n - 1, n.
\eeq

According to \eqnref{fj} and \eqnref{anomal1}, we have
\beq\label{appF2}
\sum_{j=1}^\infty \frac{\Ga_j(z)^2}{\Gd^2 + \Gl_j^2} \approx \sum_{n = 1}^{\infty} \frac{1}{\Gd^2 + \Gl_n^2} \sum_{m=-n}^n \big| a \cdot \nabla \Scal_{\p\GO}[\psi_n^m] (z) \big|^2.
\eeq
We see from \eqnref{singleball} that $\Scal_{\p\GO}[\psi_n^m]$ is a homogeneous harmonic polynomial of degree $-(n+1)$ in $\Rbb^3 \setminus \overline{\GO}$. So, we have
\beq\label{xjderi}
\pd{}{x_j} \Scal_{\p\GO}[\psi_n^m] (x)= r^{-(n+2)} \sum_{m=-n-1}^{n+1} a_{nm}^j Y_{n+1}^m (\hat{x}), \quad |x| >1
\eeq
for some constants $a_{nm}^j$, $j=1,2,3$. We then have
$$
\int_{\Rbb^3 \setminus\GO} \big| \nabla \Scal_{\p\GO}[\psi_n^m] (x) \big|^2 dx = \frac{1}{2(n+2)} \sum_{j=1}^3 \sum_{m=-n-1}^{n+1} |a_{nm}^j|^2.
$$
We then infer from \eqnref{gradest} that
$$
\frac{1}{2(n+2)} \sum_{j=1}^3 \sum_{m=-n-1}^{n+1} |a_{nm}^j|^2 \lesssim \| \psi_n^m \|_{\Hcal^*} = 1.
$$
It then follows from \eqnref{xjderi} that
\begin{align*}
\big|\nabla \Scal_{\p\GO}[\psi_n^m] (z)\big|^2 &= r_0^{-2(n+2)} \Big( \sum_{j=1}^3 \sum_{m=-n-1}^{n+1} |a_{nm}^j|^2 \Big) \Big( \sum_{m=-n-1}^{n+1} |Y_n^m (\hat{z})|^2 \Big) \\
& \lesssim (n+1) r_0^{-2(n+2)} \sum_{m=-n-1}^{n+1} |Y_n^m (\hat{z})|^2,
\end{align*}
where $r_0=|z|$. Using the Uns\"{o}ld's theorem
  \[
  \sum_{m = -n}^n \left\vert Y_n^m (\hat{z}) \right\vert^2 = \frac{2 n + 1}{4 \pi}, \ n = 0, 1, 2, \cdots,
  \]
we obtain
\beq
\big|\nabla \Scal_{\p\GO}[\psi_n^m] (z)\big|^2 \lesssim n^2 r_0^{-2(n+2)}.
\eeq

We now see from \eqnref{appF2} that
$$
\| \nabla (u_\Gd -F_z) \|_{L^2(\GO)}^2 \lesssim \sum_{n = 1}^{\infty} \frac{n^3 r_0^{-2(n+2)}}{\Gd^2 + \Gl_n^2} \lesssim \sum_{n = 1}^{\infty} n^5 r_0^{-2(n+2)} \lesssim 1
$$
for any $\Gd$. Since $\| \nabla F_z \|_{L^2(\GO)}^2$ is finite, the proof is completed. \qed

\section*{Appendix: Expansions of the fundamental solution}

In this appendix we write down the formula \eqnref{decomp} explicitly when $\GO$ is a unit ball in $\Rbb^3$ or an ellipse in $\Rbb^2$ to demonstrate the connection of \eqnref{decomp} with the known expansion formula of the fundamental solution for the Laplace equation.

Let $\GO$ be the unit ball in $\Rbb^3$. We see from \eqnref{eigenfunball} that the eigenfunctions in $\Hcal^*_0$ are
\[
\{\frac{1}{\sqrt{2n + 1}} Y_n^m ; \ n = 1, 2, \ldots\text{ and } m = - n, -n + 1, \ldots, n - 1, n,\} .
\]
One can easily see that $\Gvf_0 (x) = \frac{1}{4\pi}$ and
\beq\label{singleball_psi0}
\Scal_{\p\GO} [\varphi_0] (x) =
\begin{cases}
\ds -  \frac{1}{4\pi}, \quad & \text{ for } |x| = r < 1, \\
\nm
\ds -  \frac{1}{4\pi} \frac{1}{r}, \quad & \text{ for } |x| = r \ge 1.
\end{cases}
\eeq
Then, using \eqnref{singleball} and \eqnref{singleball_psi0}, we have
\beq\label{expansion_ball}
\GG(x-z) = - \sum_{n = 1}^{\infty} \frac{r^n r_z^{- n - 1}}{(2 n + 1)^3} \sum_{m = -n}^n Y_{n,m}(\hat{x}) Y_{n,m}(\hat{z}) - \frac{1}{4\pi} \frac{1}{r_z}, \quad r \le 1 < r_z,
\eeq
where $r = |x|$, $r_z = |z|$ and $\hat{x} = x / |x|$, $\hat{z} = z / |z|$. As mentioned right after Theorem \ref{thm:fund}, \eqref{expansion_ball} is the well-known expansion of the fundamental solution in $\Rbb^3$ in terms of spherical harmonics.

Now let $\GO$ be an ellipse with the elliptic radius $\rho_0$ in $\Rbb^2$. Then eigenfunctions in $\Hcal^*_0$ are given by $\{\psi_n^c, \psi_n^s\}_{n=1}^{\infty}$ defined by \eqnref{515}. We can see from \eqnref{415} and \eqnref{416} that $\Xi (\rho_0, \Go)^{-1}$ is orthogonal to  $\{\psi_n^c, \psi_n^s\}_{n=1}^{\infty}$ in $\Hcal^*$. So we have from \eqnref{one} that
\[
\varphi_0 (\Go) = \frac{1}{2\pi \Xi (\rho_0, \Go)}.
\]
Moreover, one can see that
\begin{align}\label{519}
  \Scal_{\p\GO} [\varphi_0] (x)
  =
  \begin{cases}
    \ds \frac{1}{2\pi} (\rho_0 + \ln {R} - \ln {2}), \quad & \rho < \rho_0, \\
    \nm
    \ds \frac{1}{2\pi} (\rho + \ln {R} - \ln {2}), \quad & \rho \ge \rho_0,
  \end{cases}
\end{align}
where $(\rho, \omega)$ be the elliptic coordinates of $x \in \Rbb^2$. In fact, as $\rho \to \infty$, we have
\[
|x| = R \sqrt{\cosh^2 {\rho} \cos^2 {\omega} + \sinh^2 {\rho} \sin^2 {\omega}} = \frac{R}{2} e^{\rho} + O(1),
\]
which implies that
\[
\rho + \ln R - \ln 2 = \ln |x| + O(|x|^{-1}) \quad\mbox{as } |x| \to \infty.
\]
Note that $\rho$ is a harmonic function. So, the righthand side of \eqnref{519} is constant in $\rho < \rho_0$, harmonic in $\rho > \rho_0$, and behaves like $\ln |x|$ as $|x| \to \infty$. So the equality in \eqnref{519} holds.
Using \eqnref{415}, \eqnref{416} and \eqnref{519}, we now have
\begin{align}\label{expansion_ellipse}
  \begin{split}
    \GG(x-z) = & - \sum_{n = 1}^{\infty} \frac{1}{n\pi} \left( \cosh{n\rho}\cos{n\rho} \, e^{- n \rho_z}\cos{n\rho_z}  + \sinh{n\rho} \sin{n\rho} \, e^{- n \rho_z}  \sin{n\rho_z}  \right) \\
    & \quad + \frac{1}{2\pi} \left(\rho_z + \ln {\left(\frac{R}{2}\right)}\right), \quad \rho \le \rho_0 < \rho_z,
  \end{split}
\end{align}
where $(\rho, \omega)$ is the elliptic coordinate of $x \in \overline{\GO}$ and $(\rho_z, \omega_z)$ is that of $z \in \Rbb^2 \setminus \overline{\GO}$.
An expansion of a derivative of $\GG(x-z)$ is also obtained in \cite[(4.6)]{CKKL}.


\end{document}